\documentclass{amsart}
\usepackage{tikz,tkz-euclide}
\title{Singularities and clusters.}
\author{V.V.Fock}
\address{Universit\'e de Strasbourg}
\begin{document}
\begin{abstract}
The aim of this note is to describe a geometric relation between simple plane curve singularities classified by simply laced Cartan matrices and cluster varieties of finite type also classified by the simply laced Cartan matrices. We construct certain varieties of configurations of flags out of Dynkin diagrams and out of singularities and show that they coincide if the Dynkin diagram corresponds to the singularity.
\end{abstract}
\dedicatory{Dedicated to Norbert A'Campo on occasion of his 80-th birthday.}
\maketitle

The correspondence between singularities\index{singularity} and cluster varieties\index{cluster variety} was first observed by S.Fomin, P.Pylyavsky, E.Shustin, and  D.Thurston in their remarkable paper \cite{FPST}.  Starting from certain real forms of deformations of the singularities introduced by N.A'Campo \cite{AC} and S.Gussein-Zade \cite{G-Z} they produced a cluster variety and showed that different resolutions of the same singularity give the same cluster variety. In this note we try to make this correspondence more geometrical and less mysterious. In particular we show that there is a map from the base of a versal deformation of the singularity to the corresponding cluster variety. For this purpose we generalize a construction by R.Nevanlinna \cite{N} brought to our attention by B.Shapiro \cite{S}. R.Nevanlinna studied the map from the space of the differential operators of the form $D=\partial^2/\partial z^2+P(x)$, where $P$ is a polynomial, to the collections of points in a projective line $P^1$. These points just correspond to solutions of the equation $D\psi=0$ vanishing at infinity along different rays on the complex plane and viewed as lines in the two-dimensional space of global solutions. On the other hand, the symbol of this operator $p^2+P(x)$ is just an equation for the versal deformation of a plane curve singularity of type $A$. 

Our construction is in a sense a generalization of this one for singularities of other types. Namely, a versal family of a planar singularity can be considered as a family of Lagrangian curves in the complex symplectic plane $(\mathbb{C}^2,dp\,dx)$. Observe that this family can be locally parameterized by the cohomology class of the form $pdx$. On the other hand the family of equations defining these curves can be transformed into a family of differential operators having the equations as their symbol. The space of Stokes data at infinity of the differential operator coincides with the cluster variety suggested in \cite{FPST}. 

The correspondence between symbols and operators is of course not canonical. First of all, we can change a representation of the symplectic plane $(\mathbb{C}^2,dp\,dx)$ as a cotangent bundle to a complex line $\mathbb{C}$ in different ways. All these ways give different Stokes data, but, as one can verify in the examples, they are equivalent as cluster varieties. Sometimes the birational equivalence of the corresponding configuration spaces is not so trivial and can be considered as a generalization of the Gale duality.

On the other hand the correspondence between  the space of differential operators and its symbol is not canonical either. We conjecture that this map becomes canonical in the tropical limit and the cluster coordinates in this limit are periods of the Lagrangian curves given by the versal family of the singular one. 

I am strongly indebted to F.Chapoton and E.Shustin for reading the paper and making very crucial remarks.

\section{Recall: Configurations corresponding to planar bipartite graphs}

Let $a$ be a collection of points of a projective space $P(V)$ vector space $V$. Denote by $\langle a\rangle$ a projective subspace generated by the points from $a$. Recall that a collection $a$ of $k$ points is called \textit{free} if $\dim\, \langle a\rangle=k-1$ and a \textit{circuit} if $\dim\, \langle a\rangle=k-2$. For example, two points form a circuit\index{circuit} if they coincide, three points form a circuit if they are colinear four points form a circuit if they are coplanar \textit{etc.}

Let $\Gamma$ be a bipartite graph with the set of white vertices $W$, black vertices $B$ and edges $E$. For simplicity we assume here that any two vertices are connected by no more than one edge. We say that an association of a point $p_w$ in a projective space to every white vertex $w\in W$ \textit{corresponds} to $\Gamma$ if for every black vertex $b\in B$ the points corresponding to its white neighbors form a circuit. We also require that the collection of points does not belong to a proper projective subspace.

Denote by $\mathcal{X}_\Gamma$ the set of configurations corresponding to $\Gamma$ considered up to the action of the projective group $PGL(V)$. We call the \textit{dimension of the configuration} the dimension of the projective space $P(V)$. We say that the graph $\Gamma$ is \textit{minimal} if removal of any black vertex increases the dimension of the corresponding configuration.

A discrete connection on a graph $\Gamma$ is an association of one dimensional vector spaces to vertices and association to edges of isomorphisms of the spaces corresponding to its endpoints. Given a basis in each of the one-dimensional subspaces, the discrete connection becomes an association of nonzero numbers to edges. These numbers can be organized into a matrix $M_b^w$ with columns and rows enumerated by the black and white vertices, respectively and with zeroes for pairs of vertices not connected by an edge. Changing the bases amounts to the multiplication of this matrix by invertible diagonal matrices from the left and from the right. Given a closed path of the graph, the monodromy of the connection around this graph is a composition of the maps corresponding to the edges. In terms of the connection matrix, if the path passes consecutively through the vertices $b_0,w_0,\ldots,b_k,w_k,b_0$ then the monodromy is given by a Laurent monomial in the matrix entries $M_{b_1}^{w_1}(M_{b_2}^{w_1})^{-1}\cdots(M_{b_0}^{w_k})^{-1}$. If the graph is planar a graph connection is uniquely determined by the monodromies around its faces.

The set of graph connections  on a graph $\Gamma$ is in a bijection with the set of configurations corresponding to $\Gamma$. Indeed, given the matrix $M_b^w$ representing the connection we can consider it as a map $M:\mathbb{C}^B\to \mathbb{C}^W$. The image of the standard basis of $\mathbb{C}^W$ in the projectivized cokernel $P(\mathbb{C}^W/ \textrm{Im }M)$ form the desired configuration. Conversely, given a configuration of points in a projective space $P(V)$, for every white vertex $w$ choose a vector $\tilde{p}_w\in V$ representing $p_w$. Given a black vertex $b$ the chosen vectors corresponding to its white neighbors satisfy a nontrivial linear relation $\sum M_b^w \tilde{p}_w=0$. The matrix $(M_b^w)$ forms the desired graph connection. 

Call two configurations \textit{equivalent}\index{equivalent configurations} if one can construct one from another by adding a point which can be constructed from the remaining points by tracing subspaces through given points and taking their intersections. We can also remove a point if its position is uniquely determined by other points. The graphs corresponding to equivalent configurations are also called equivalent. Many examples of such equivalent configurations will be given below. 

There are four equivalences generating many  (and conjecturally all) others. Namely 
\begin{enumerate}
\item[1.]
\raisebox{-.2\height}{
\begin{tikzpicture}
\draw (-0.5,0)--(0.5,0);
\draw (-0.7,0.2)--(-0.5,0)--(-0.7,-0.2);
\draw (0.7,0.2)--(0.5,0)--(0.7,-0.2);
\fill (0,0) circle (2pt);
\draw[fill=white] (-0.5,0) circle (2pt);
\draw (-0.75,0.1) node {$\vdots$};
\draw (0.75,0.1) node {$\vdots$};
\draw[fill=white] (0.5,0) circle (2pt);
\draw (1.25,0) node {$\longleftrightarrow$};
\draw (1.8,0.2)--(2.2,-0.2);
\draw (1.8,-0.2)--(2.2,0.2);
\draw[fill=white] (2,0) circle (2pt);
\draw (2.25,0.1) node {$\vdots$};
\draw (1.75,0.1) node {$\vdots$};
\end{tikzpicture}
}

\noindent Two-valent black vertex means that two points corresponding to its white neighbors coincide, and of course one of them can be removed. 
\item[$1'$.]
\raisebox{-.2\height}{
\begin{tikzpicture}
\draw (-0.5,0)--(0.5,0);
\draw (-0.7,0.2)--(-0.5,0)--(-0.7,-0.2);
\draw (0.7,0.2)--(0.5,0)--(0.7,-0.2);
\draw[fill=white] (0,0) circle (2pt);
\fill (-0.5,0) circle (2pt);
\draw (-0.75,0.1) node {$\vdots$};
\draw (0.75,0.1) node {$\vdots$};
\fill (0.5,0) circle (2pt);
\draw (1.25,0) node {$\longleftrightarrow$};
\draw (1.8,0.2)--(2.2,-0.2);
\draw (1.8,-0.2)--(2.2,0.2);
\fill (2,0) circle (2pt);
\draw (2.25,0.1) node {$\vdots$};
\draw (1.75,0.1) node {$\vdots$};
\end{tikzpicture}
}

\noindent Two-valent white vertex means that the circuits corresponding to its black neighbors have one common point, and obviously their union (with one of the points de-doubled) is also a circuit.
\item[2.]
\raisebox{-.45\height}
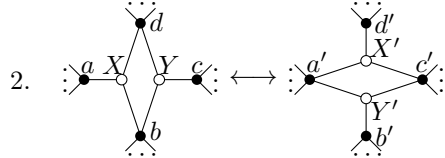
{
\begin{tikzpicture}
\draw (-0.75,0)--(-0.25,0);
\draw (0.75,0)--(0.25,0);
\draw (0.25,0)--(0,0.75)--(-0.25,0)--(0,-0.75)--cycle;
\draw (0.95,0.2)--(0.75,0)--(0.95,-0.2);
\draw (-0.95,0.2)--(-0.75,0)--(-0.95,-0.2);
\draw (-0.2,0.95)--(0,0.75)--(0.2,0.95);
\draw (-0.2,-0.95)--(0,-0.75)--(0.2,-0.95);
\draw (-1,0.1) node {$\vdots$};
\draw (1,0.1) node {$\vdots$};
\draw (0.05,1) node {$\cdots$};
\draw (0.05,-1) node {$\cdots$};
\fill (0,0.75) circle (2pt);\draw (0.21,0.77) node{$d$};
\fill (0,-0.75) circle (2pt);\draw (0.2,-0.68) node{$b$};
\fill (0.75,0) circle (2pt);\draw (0.75,0.17) node{$c$};
\fill (-0.75,0) circle (2pt);\draw (-0.70,0.17) node{$a$};
\draw[fill=white] (-0.25,0) circle (2pt);\draw (-0.35,0.2) node{$X$};
\draw[fill=white] (0.25,0) circle (2pt);\draw (0.38,0.2) node{$Y$};
\draw (1.5,0) node {$\longleftrightarrow$};
\draw (3,-0.75)--(3,-0.25);
\draw (3,0.75)--(3,0.25);
\draw (3,0.25)--(3.75,0)--(3,-0.25)--(2.25,0)--cycle;
\draw (3.95,0.2)--(3.75,0)--(3.95,-0.2);
\draw (2.05,0.2)--(2.25,0)--(2.05,-0.2);
\draw (2.8,0.95)--(3,0.75)--(3.2,0.95);
\draw (2.8,-0.95)--(3,-0.75)--(3.2,-0.95);
\draw (2,0.1) node {$\vdots$};
\draw (4,0.1) node {$\vdots$};
\draw (3.05,1) node {$\cdots$};
\draw (3.05,-1) node {$\cdots$};
\fill (3,0.75) circle (2pt);\draw (3.25,0.75) node{$d'$};
\fill (3,-0.75) circle (2pt);\draw (3.25,-0.75) node{$b'$};
\fill (3.75,0) circle (2pt);\draw (3.78,0.24) node{$c'$};
\fill (2.25,0) circle (2pt);\draw (2.35,0.24) node{$a'$};
\draw[fill=white] (3,-0.25) circle (2pt); \draw (3.25,-0.4) node{$Y'$};
\draw[fill=white] (3,0.25) circle (2pt); \draw (3.25,0.43) node{$X'$};
\end{tikzpicture}
}

\noindent Let $a$,$b$,$c$,$d$ be the circuits corresponding to black vertices and $X$ and $Y$ points corresponding to the white ones.  Let $Y'=\langle b\setminus \{X,Y\}\rangle\cap (XY)$. The intersection is a point since obviously the line $(XY)$ belongs to the subspace $\langle b\rangle$ and the subspace $\langle b\setminus \{X,Y\}\rangle$ is of codimension 1 in it. Similarly let $X'=\langle d\setminus \{X,Y\}\rangle\cap (XY)$. 
We can now define $a'=a\cup \{X',Y'\}\setminus\{X\}$,  $b'=b\cup \{Y'\}\setminus\{X,Y\}$, $c'=c\cup \{X',Y'\}\setminus\{Y\}$, $d'=d\cup \{X'\}\setminus\{X,Y\}$.
\item[$2'$.]
\raisebox{-.45\height}
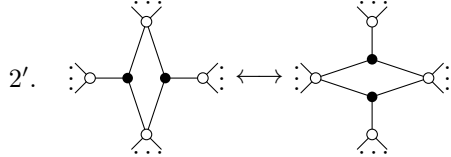
{
	\begin{tikzpicture}
	\draw (-0.75,0)--(-0.25,0);
	\draw (0.75,0)--(0.25,0);
	\draw (0.25,0)--(0,0.75)--(-0.25,0)--(0,-0.75)--cycle;
	\draw (0.95,0.2)--(0.75,0)--(0.95,-0.2);
	\draw (-0.95,0.2)--(-0.75,0)--(-0.95,-0.2);
	\draw (-0.2,0.95)--(0,0.75)--(0.2,0.95);
	\draw (-0.2,-0.95)--(0,-0.75)--(0.2,-0.95);
	\draw (-1,0.1) node {$\vdots$};
	\draw (1,0.1) node {$\vdots$};
	\draw (0.05,1) node {$\cdots$};
	\draw (0.05,-1) node {$\cdots$};
	\draw[fill=white] (0,0.75) circle (2pt);
	\draw[fill=white] (0,-0.75) circle (2pt);
	\draw[fill=white] (0.75,0) circle (2pt);
	\draw[fill=white] (-0.75,0) circle (2pt);
	\fill (-0.25,0) circle (2pt);
	\fill (0.25,0) circle (2pt);
	\draw (1.5,0) node {$\longleftrightarrow$};
	\draw (3,-0.75)--(3,-0.25);
	\draw (3,0.75)--(3,0.25);
	\draw (3,0.25)--(3.75,0)--(3,-0.25)--(2.25,0)--cycle;
	\draw (3.95,0.2)--(3.75,0)--(3.95,-0.2);
	\draw (2.05,0.2)--(2.25,0)--(2.05,-0.2);
	\draw (2.8,0.95)--(3,0.75)--(3.2,0.95);
	\draw (2.8,-0.95)--(3,-0.75)--(3.2,-0.95);
	\draw (2,0.1) node {$\vdots$};
	\draw (4,0.1) node {$\vdots$};
	\draw (3.05,1) node {$\cdots$};
	\draw (3.05,-1) node {$\cdots$};
	\draw[fill=white] (3,0.75) circle (2pt);
	\draw[fill=white] (3,-0.75) circle (2pt);
	\draw[fill=white] (3.75,0) circle (2pt);
	\draw[fill=white] (2.25,0) circle (2pt);
	\fill (3,-0.25) circle (2pt);
	\fill (3,0.25) circle (2pt);
	\end{tikzpicture}
}

\noindent This fragment of the graph corresponds to four collinear points. Any triple of such point is a circuit. Any two of such circuits define such configurations and can be exchanged by the remaining two preserving planarity of the graph.
\end{enumerate}

\section{From Dynkin diagram to bipartite graph.} For any planar graph (in our examples they will be Dynkin diagrams) one can associate a planar bipartite graph dual to the original graph (See example shown on Fig. \ref{an}a, \ref{d4}a and \ref{e8}a.) Namely, surround every vertex of valency $k>2$ by a $2k$-gon in a way that its every second side crosses the respective edge of the original graph. For a 1-valent vertex surround it with a square. Then identify the sides of the polygons crossing the same edges of the original graph. Finally color vertices in black and white and add one edge with a white vertex to the every black 2-valent vertex drawn around vertices of valency 1. 

\section{Stokes data.}
Let $D(z,\frac{\partial}{\partial z})$ be a differential operator with polynomial coefficients of order $n$. Let $D(x,p)$ be its symbol and $\Delta$ be its Newton polygon\index{Newton polygon}, namely the convex hull of the points with coordinates $(a,b)\in \mathbb{Z}^2$ for each monomial term $x^ap^b$ of the polynomial $D(x,p)$. 

Consider the algebraic curve $C=\{(x,p)\in (\mathbb{C}^*)^2| D(x,p)=0\}$.  This curve has the genus equal to the number of integer points strictly inside the polygon $\Delta$ and can be compactified by adding points, called \textit{compactification points}, corresponding to sides of the polygon. By a side we mean a segment of the boundary of $\Delta$ between two adjacent points of $\mathbb{Z}^2$.  If we orient the boundary of the polygon counterclockwise each side $s$ correspond to indivisible vector $(a,b)\in \mathbb{Z}^2$. The sum of such vectors obviously vanishes. At the compactification points corresponding to $s$ the functions $x$ and $p$ have a zero of order $-b_s$ and $a_s$, respectively. Therefore in the vicinity of such point we have $x^{a_s}p^{b_s}=O(1)$ and hence $p\sim A_sx^{-a_s/b_s}$ for some constants $A_s$.

One should however realize that the correspondence between sides and compactification points is not canonical but is defined up to permutations of sides corresponding to equal vectors.

The same polygon defines a curve in $\mathbb{C}^2$. The sides of the polygon with $a>0$ and $b\leq 0$ correspond to compactifcation points with both coordinates tending to a finite constant and thus belonging to the curve in $\mathbb{C}^2$. The remaining compactification points are called \textit{points at infinity}.

Consider now the equation $D(x,p)=0$ as a family of equations for the indeterminate $p$ depending on $x$ as a parameter. Our aim is to determine the asymptotic behavior of its roots when $x\to \infty$. This behavior of $x$ corresponds to the sides of the polygon with $b_s>0$ (and thus going upward on a picture where the $p$-axis points up). Each side corresponds to $b_s$ roots with asymptotic $p\sim A_s x^{-a_s/b_s}$. The total number of roots is thus equal to the height of the polygon, which is just the degree of the polynomial $D(x,p)$ with respect to $p$. 

Now proceed to the study of the asymptotic behaviour of the differential equation $D(x,\frac{\partial}{\partial x})\psi =0$. In a simply connected domain where the coefficient of the $n$-th derivative vanishes nowhere, it has an $n$-dimensional space of solutions. However in a sufficiently small vicinity of infinity the differential equation defines an $n$-dimensional local system. Indeed, fix a number $R$ such that the coefficient of the highest derivative of $D$ does not vanish for $|x|>R$, a number $\varepsilon\in ]0,\pi/2[$ and an angle $\alpha\in \mathbb{R}/2\pi\mathbb{Z}$. The domain $\{x\in \mathbb{C}\vert\ |\arg x-\alpha|<\varepsilon,\ |x| > R\}$ satisfies the above conditions. Define $V_\alpha$ as the space of solutions of the differential equation in this domain. The space $V_\alpha$ does depend neither of $\varepsilon$ nor of $R$ and does not depend on $\alpha$ locally. Namely if $|\alpha_1-\alpha_2|<\varepsilon$ the corresponding domains intersect, their intersection is simply connected and therefore the spaces $V_{\alpha_1}$ and $V_{\alpha_2}$ can be identified. Hence the family $V_\alpha$ form an $n$-dimensional local system over the circle $\mathbb{R}/2\pi\mathbb{Z}$.

Every solution of the equation has asymptotic behavior $\psi(x) \sim e^{\int{p(x)dx}}$ for a root $p(x)$ and therefore for generic angle $\alpha$ the space $V_\alpha$ is filtered by the rates of growth of the functions $\Re e\int^{re^{i\alpha}}\!\!p(x)dx$ with $r\to +\infty$. 

Consider a sufficiently large $R$ and mark the points $A_s(Re^{i\alpha})^{(b_s-a_s)/b_s}$ on the complex plane. (For a given $s$ there are $b_s$ such points.) These points are ordered according to their projection on the real axis and the interval between the $i$-th and the $i+1$-st points correspond to the subspace of $V_\alpha$ of dimension $i$. As $\alpha$ runs around the circle, the points rotate around the origin with the angular speed $(b_s-a_s)/b_s$ and the order of their projections changes. We call the collection of points corresponding to generic $\alpha$ with their angular speeds indicated for each point the \textit{growth diagram}.  When the $i$-th and the $i+1$-th projections pass through each other the $i$-th dimensional subspace of $V_\alpha$ changes. 

The local system and the sequence of flags in its fibers constitute the Stokes data of the differential operator $D$ at infinity. If the coefficient of the highest derivative in the operator $D$ is constant, the local system is trivial and the Stokes data at infinity amounts to the collection of flags in a fixed vector space. 

For nontrivial local systems, taking the universal cover of the circle we can consider a finite sequence of flags in a local system as an infinite quasi-periodic one in a fixed vector space. Recall that the sequence is called quasi-periodic if its shift by a period coincides with the action of an element of $GL(n)$. 

Recall that pairs of complete flags in an $n$-dimensional space up to diagonal action of $GL(n)$ are in bijection with the permutation group $\mathfrak{S}_n$. The standard generators $s_i$ of this group correspond to the pairs of flags different only in the subspaces of dimension $i$. 

Therefore a sequence of flags in a local system on a circle such that adjacent flags differ in one subspace can be encoded by an infinite periodic word in the same generators $s_1,\ldots,s_{n-1}$. We will denote such words as $[w]$, where $w$ is a period of the infinite word. For example a word $[(s_2s_1)^m]=[(s_1s_2)^m]$ corresponds to a 3-dimensional local system with $2m$ flags such that the subspaces of dimension $1$ and $2$ change alternatively. Such sequence is equivalent to the quasiperiodic sequence of 1-dimensional subspaces since the 2-dimensional subspaces can be restored from the 1-dimensional ones. In the projective space this sequence corresponds to a quasi-periodic set of points with the period of length $m$.

Observe that if two words are related by braid relations they correspond to equivalent sequences of flags. 

The word $[(s_2s_1^2)^m]$ corresponds to a sequence of flags where one dimensional subspace changes twice after each change of the 2-dimensional subspace. In the projective space it corresponds to a quasi-periodic broken line with a marked point on each side.

As another example consider a growth diagram given by a regular $n$-gon rotating around its center with angular velocity $m/n$.
We will show that it corresponds to an $m$-periodic sequence of points in $P^{n-1}$. Indeed, in this case all even subspaces change at once and then all odd subspaces change at once. The corresponding word is therefore $w=[(w_\textrm{odd}w_\textrm{even})^m]$ 
where $w_{\textrm{odd}}=s_1s_3\cdots$ and $w_{\textrm{even}}=s_2s_4\cdots$ are the products of the odd and even generators, respectively. We claim that such configuration is equivalent just to the sequence of 1-dimensional subspaces with quasi-period $m$. Indeed, given a sequence of flags corresponding to the word $w$ one can construct the sequence of 1-dimensional subspaces just discarding all subspaces of higher dimension. On the other hand given a sequence of 1-dimensional subspaces $\{\tilde{p}_i|i\in \mathbb{Z}/n\mathbb{Z}\}$, one can construct a sequence of $2m$ flags 

$$F_{2i+1}=\{ \tilde{p}_i\subset \tilde{p}_i+\tilde{p}_{i+1}\subset  \tilde{p}_{i-1}+\tilde{p}_i+\tilde{p}_{i+1}\subset  \tilde{p}_{i-1}+\tilde{p}_i+\tilde{p}_{i+1}+\tilde{p}_{i+2}\subset \cdots\} $$
and
$$F_{2i}=\{ \tilde{p}_{i}\subset \tilde{p}_{i-1}+\tilde{p}_{i}\subset  \tilde{p}_{i-1}+\tilde{p}_i+\tilde{p}_{i+1}\subset  \tilde{p}_{i-2}+\tilde{p}_{i-1}+\tilde{p}_{i}+\tilde{p}_{i+1}\subset \cdots\} $$
and observe that $F_{2i}$ differ from $F_{2i+1}$ in odd dimensional terms and from $F_{2i-1}$ in even dimensional terms.

\section{Example: $A_n$.}
The Dynkin diagram of type $A_n$ is just a chain of $n$ vertices and the corresponding bipartite graph shown on Fig.\ref{an}a consists of $n+3$ white and $n+1$ black vertices. 
\begin{figure}[h]
\begin{tikzpicture}
\draw[thick,green] (0,0)--(1,0);
\draw[thick,green,dashed] (1,0)--(2.5,0);
\fill[green] (0,0) circle(2pt);
\fill[green] (0.5,0) circle(2pt);
\fill[green] (1,0) circle(2pt);
\fill[green] (2.5,0) circle(2pt);
\fill[green] (3,0) circle(2pt);
\draw[green] (2.5,0)--(3,0);
\draw (-0.25,-0.25)--(-0.25,0.25)--(0.25,0.25)--(0.25,-0.25)--(-0.25,-0.25)--(-0.6,-0.6);
\draw (0.25,0.25)--(0.75,0.25)--(0.75,-0.25)--(0.25,-0.25);
\draw (0.75,0.25)--(1.25,0.25)--(1.25,-0.25)--(0.75,-0.25);
\draw[dashed] (1.25,0.25)--(2.25,0.25);
\draw[dashed] (1.25,-0.25)--(2.25,-0.25);
\draw (2.25,0.25)--(2.75,0.25)--(2.75,-0.25)--(2.25,-0.25)--cycle;
\draw (2.75,0.25)--(3.25,0.25)--(3.25,-0.25)--(2.75,-0.25);
\draw (3.6,0.6)--(3.25,0.25);
\foreach \x in {1,3,5,7} \fill (\x/2-0.25,0.25) circle(2pt);
\foreach \x in {0,2,6} \fill (\x/2-0.25,-0.25) circle(2pt);
\foreach \x in {0,2,6} \draw[fill=white] (\x/2-0.25,0.25) circle(2pt);
\foreach \x in {1,3,5,7} \draw[fill=white] (\x/2-0.25,-0.25) circle(2pt);
\draw[fill=white] (-0.6,-0.6) circle(2pt);
\draw[fill=white] (3.6,0.6) circle(2pt);
\draw (1.5,-1.2) node {a};
\end{tikzpicture}\quad
\begin{tikzpicture}
\draw[->] (-0.25,0) -- (2.75,0);
\draw[->] (0,-0.25) -- (0,1.25);
\fill[color=gray!20](0,1)--(0,0)--(2.5,0)--cycle;
\draw[thick] (0,1)--(0,0)--(2.5,0);
\draw[thick,blue](0,1)--(2.5,0);
\foreach \x in {0,1,2,3,4,5} \fill (\x/2,0) circle(2pt);
\foreach \x in {0,1,2} \fill (\x/2,0.5) circle(2pt);
\foreach \x in {0} \fill (\x/2,1) circle(2pt);
\fill (0,1) circle(1.5pt);
\fill[white] (0.5,1) circle(1.5pt);
\fill[white] (2,0) circle(1.5pt);
\draw (1.25,-0.8) node {b};
\end{tikzpicture}\quad
\begin{tikzpicture}
\draw[->] (-0.9,0)--(0.9,0);
\draw[->] (0,-0.9)--(0,0.9);
\draw[-latex,blue] (-70:0.8)arc(-70:-50:0.8);
\draw[-latex,blue] (110:0.8)arc(110:130:0.8);
%\draw[-latex,blue] (170:0.8)arc(170:190:0.8);
\fill[blue] (-70:0.8) circle(1.5pt);
%\fill[blue]  (50:0.8) circle(1.5pt);
\fill[blue]  (110:0.8) circle(1.5pt);
\draw (-45:0.8) node {\tiny $(n+3)/2$};
\draw (137:0.8) node {\tiny $(n+3)/2$};
\draw (0,-1.4) node {c};
\end{tikzpicture}
\caption{}\label{an}
\end{figure}
It corresponds just to $n+3$ collinear points. 

The versal deformation of a singularity of type $A_n$ is represented by the polynomial
$$ D(x,p)=p^2+x^{n+1}+a_1+a_2x+\cdots+a_nx^{n-1}.
$$
The corresponding Newton polygon is a right triangle with sides $2$ and $n$ (shown on Fig.\ref{an}b for $n=4$) with one or two sides for $n$ odd or even, respectively, going upwards. It corresponds to a curve of genus $\lfloor (n-1)/2\rfloor$ with one or two points at infinity, respectively with homology of rank $n$. The growth diagram (shown on Fig.\ref{an}c) consists of two points, which are opposite if $n$ is even, rotating with the angular speed $(n+3)/2$. It corresponds to the word $[s_1^{n+3}]$ and thus the configurations of $n+3$ points in $P^1$.

Interchanging $p$ and $x$ we get a differential equation of order $n+1$, the growth diagram consists of points forming a regular $n+1$-gon rotating about its center with angular velocity $(n+3)/(n+1)$ which corresponds to the cyclic word $[ (w_\textrm{odd}w_\textrm{even})^{n+3}]$ and thus to $n+3$ points in $P^{n}$. These two configuration space are known to be isomorphic since one can trace a unique rational normal curve through this point thus obtaining $(n+3)$ points in $P^1$.

\section{Example: $D_4$.}
Consider the Dynkin diagram of the type $D_4$ and construct a bipartite graph shown on Fig. \ref{d4}a. It corresponds to three groups of collinear points $C,1,2,A$, $A,3,4,B$ and $B,5,6,A$ and the dimension of the configuration is 2. The configuration is obviously equivalent to the configuration of 6 points $1,2,3,4,5,6 \in P^2$ as shown on Fig.\ref{d4}b. This configuration of points corresponds to a cyclic word $[(s_2s_1^3)^3]$
\begin{figure}[h]
\begin{tikzpicture}
\draw (-0.3,0)--(0,0.5)--(0.6,0.5)--(0.9,0)--(0.6,-0.5)--(0,-0.5)--cycle;
\draw (0,-0.5)--(0,-1.1)--(0.6,-1.1)--(0.6,-0.5);
\draw (0.6,0.5)--(1.1,0.8)--(1.4,0.3)--(0.9,0);
\draw (-0.3,0)--(-0.8,0.3)--(-0.5,0.8)--(0,0.5);
\draw (1.1,0.8)--(1.31,1.36);
\draw (-1.37,0.16)--(-0.8,0.3);
\draw (0.6,-1.1)--(1.02,-1.52);
\draw[green] (-0.4,0.4)--(0.3,0)--(1,0.4);
\draw[green] (0.3,0)--(0.3,-0.8);
\fill[green] (-0.4,0.4) circle (2pt);
\fill[green] (0.3,0) circle (2pt);
\fill[green] (0.3,-0.8) circle (2pt);
\fill[green] (1,0.4) circle (2pt);
\fill (0,0.5) circle (2pt);
\fill (-0.8,0.3) circle (2pt);
\fill (0,-0.5) circle (2pt);
\fill (0.6,-1.1) circle (2pt);
\fill (0.9,0) circle (2pt);
\fill (1.1,0.8) circle (2pt);
\draw[fill=white] (1.02,-1.52) circle (2pt) node[below right]{$6$};
\draw[fill=white] (1.31,1.36) circle (2pt) node[above]{$2$};
\draw[fill=white] (-1.37,0.16) circle (2pt) node[left]{$4$};
\draw[fill=white] (-0.3,0) circle (2pt) node[below left]{$B$};
\draw[fill=white] (0.6,0.5) circle (2pt) node[above]{$A$};
\draw[fill=white] (0.6,-0.5) circle (2pt) node[right]{$C$};
\draw[fill=white] (1.4,0.3) circle (2pt) node[right]{$1$};
\draw[fill=white] (-0.5,0.8) circle (2pt) node[above]{$3$};
\draw[fill=white] (0,-1.1) circle (2pt) node[below left]{$5$};
\draw (0.3,-2) node {a};
\end{tikzpicture}
\begin{tikzpicture}
\tkzDefPoints{0/0/A,3/0/B,2/2/C};
\tkzDrawLine[very thick](A,B);
\tkzDrawLine[very thick](B,C);
\tkzDrawLine[very thick](C,A);
\tkzDefBarycentricPoint(B=1,C=2); \tkzGetPoint{L6};
\tkzDefBarycentricPoint(B=2,C=1); \tkzGetPoint{L5};
\tkzDefBarycentricPoint(C=2,A=4); \tkzGetPoint{L2};
\tkzDefBarycentricPoint(C=4,A=2); \tkzGetPoint{L1};
\tkzDefBarycentricPoint(A=1,B=3); \tkzGetPoint{L4};
\tkzDefBarycentricPoint(A=3,B=2); \tkzGetPoint{L3};
\tkzLabelPoint[above left=-2pt](L2){$2$};
\tkzLabelPoint[above left=-2pt](L1){$1$};
\tkzLabelPoint[below](L4){$4$};
\tkzLabelPoint[below](L3){$3$};
\tkzLabelPoint[above right=-2pt](L6){$6$};
\tkzLabelPoint[above right=-2pt](L5){$5$};
\tkzLabelPoint[below](A){$A$};
\tkzLabelPoint[below left=-2pt](B){$B$};
\tkzLabelPoint[right](C){$C$};
\tkzDrawPoints[fill = white, size=3.5pt, ](A,B,C);
\tkzDrawPoints[fill = white, size=3.5pt, color=red ](L1,L2,L3,L4,L5,L6);
\draw (1.5,-1.2) node {b};
\end{tikzpicture}
\begin{tikzpicture}[x={(0cm,1cm)},y={(1cm,0cm)}]
\draw[->] (-0.25,0) -- (1.75,0);
\draw[->] (0,-0.25) -- (0,1.25);
\fill[color=gray!20] (0,0)--(0,0.5)--(0.5,1)--(1.5,0)--cycle;
\draw[thick] (0,0)--(0,0.5)--(0.5,1)--(1.5,0)--cycle;
\draw[thick, blue] (1.5,0)--(0.5,1);
\draw[thick, magenta](0.5,1)--(0,0.5);
\foreach \x in {0,1,2,3} \fill (\x/2,0) circle(2pt);
\foreach \x in {0,2} \fill (\x/2,0.5) circle(2pt);
\foreach \x in {1} \fill (\x/2,1) circle(2pt);
\draw[fill=white] (0.5,0.5) circle(2pt);
\draw (-1.5,0.5) node {c};
\end{tikzpicture}
\begin{tikzpicture}
\coordinate (A) at (30:0.2);
\coordinate (C) at (155:0.9);
\draw[->] (-1,0)--(1,0);
\draw[->] (0,-1)--(0,1);
\draw[-latex,blue] (-60:0.8)arc(-60:-40:0.8);
\draw[-latex,blue] (145:0.9)arc(145:165:0.9);
\fill[magenta] (A) circle(1.5pt);
\fill[blue] (-60:0.8) circle(1.5pt);
\fill[blue] (145:0.9) circle(1.5pt);
\draw (-35:0.8) node {\tiny 2};
\draw (170:0.9) node {\tiny 2};
\draw (0,-2) node {d};
\end{tikzpicture}
\caption{}\label{d4}
\end{figure}

The versal deformation of the singularity $D_4$ reads as

$$D(p,x)=p^3+x^2p + a_1+a_2p+a_3p^2+a_4x$$
with the Newton polygon is a quadrilateral shown on the Fig.\ref{d4}c. There are three sides of the polygon directed upward: $(-1,1)$ with multiplicity $2$ and $(1,1)$ with multiplicity  $1$. It corresponds to a curve of genus 1 with three points at infinity, and with homology of rank 4. The growth diagram consist of two points rotating with angular speed 2 and one point closer to the center which does not move. It gives a word $[(s_2s_1^2)^4]$. Therefore the space Stokes data can be considered as the configuration space of quadrilaterals with marked points on each side. Such configurations are equivalent to the configurations of six points --- four points on the sides and two opposite vertices of the quadrilateral. On the other hand one can deduce this equivalence from the equality of the words $[(s_2s_1^2)^4]$ and $[(s_2s_1^3)^3]$ in the braid group. 

Consider another form of the same singularity just with the variables $p$ and $x$ interchanged. 

$$D(p,x)=x^3+xp^2 + a_1+a_2x+a_3x^2+a_4p$$
The Newton polygon is just the reflection of the original one, but the corresponding differential equation is of order $2$ and the local system is nontrivial since the coefficient at the highest derivative vanish at the origin. The growth diagram consists of two points rotating with the angular speed $2$ and thus the Stokes data amounts to the configuration of four lines in a two-dimensional local system on a circle.  

Remarkably these two configuration space turns out to be  birationally isomorphic. The only isomorphism I know is given by describing both as cluster varieties, and I don't know any geometric way to describe it. 

\section{Example: $E_8$.}

Consider the Dynkin diagram of the type $E_8$ and construct a bipartite graph shown on Fig. \ref{e8}a. On the same picture we show the corresponding bipartite graph. This diagram corresponds to configurations of 13 points corresponding to white vertices with black vertices corresponding to collinear triples of points. It implies that there are three groups of collinear points $C,1,2,3,4,5,A$, $A,6,7,8,B$ and $B,9,10,C$. Such configuration can be realized in two dimensional projective space $P^2$ as a triangle with 2,3 and 5 points on its sides, respectively as shown on Fig. \ref{e8}b. Observe that this configuration space is birationally isomorphic to the space of unrestricted 8-tuples of points $1,2,7,8,9,10,X,Y$ in $P_2$. Indeed , as it is clear from the picture the points $A,B,C,3,4,5,6,7$ can be reconstructed out of $1,2,7,8,9,10,X,Y$ and vise versa.

\begin{figure}[h]
\begin{tikzpicture}
\draw (-0.3,0)--(0,0.5)--(0.6,0.5)--(0.9,0)--(0.6,-0.5)--(0,-0.5)--cycle;
\draw (0,-0.5)--(0,-1.1)--(0.6,-1.1)--(0.6,-0.5);
\draw (0.6,0.5)--(1.1,0.8)--(1.4,0.3)--(0.9,0);
\draw (1.1,0.8)--(1.6,1.1)--(1.9,0.6)--(1.4,0.3);
\draw (1.6,1.1)--(2.1,1.4)--(2.4,0.9)--(1.9,0.6);
\draw (2.1,1.4)--(2.6,1.7)--(2.9,1.2)--(2.4,0.9);
\draw (-0.3,0)--(-0.8,0.3)--(-0.5,0.8)--(0,0.5);
\draw (-0.8,0.3)--(-1.3,0.6)--(-1,1.1)--(-0.5,0.8);
\draw (2.9,1.2)--(3.46,1.06);
\draw (-1,1.1)--(-1.14,1.66);
\draw (0.6,-1.1)--(1.02,-1.52);
\draw[green] (-0.9,0.7)--(0.3,0)--(2.5,1.3);
\draw[green] (0.3,0)--(0.3,-0.8);
\fill[green] (-0.9,0.7) circle (2pt);
\fill[green] (-0.4,0.4) circle (2pt);
\fill[green] (0.3,0) circle (2pt);
\fill[green] (0.3,-0.8) circle (2pt);
\fill[green] (2.5,1.3) circle (2pt);
\fill[green] (2,1) circle (2pt);
\fill[green] (1.5,0.7) circle (2pt);
\fill[green] (1,0.4) circle (2pt);
\fill (0,0.5) circle (2pt);
\fill (-0.8,0.3) circle (2pt);
\fill (-1,1.1) circle (2pt);
\fill (0,-0.5) circle (2pt);
\fill (0.6,-1.1) circle (2pt);
\fill (0.9,0) circle (2pt);
\fill (1.1,0.8) circle (2pt);
\fill (1.9,0.6) circle (2pt);
\fill (2.1,1.4) circle (2pt);
\fill (2.9,1.2) circle (2pt);
\draw[fill=white] (1.02,-1.52) circle (2pt) node[below right]{$10$};
\draw[fill=white] (-1.14,1.66) circle (2pt) node[above]{$7$};
\draw[fill=white] (3.46,1.06) circle (2pt) node[right]{$3$};
\draw[fill=white] (-0.3,0) circle (2pt) node[below left]{$B$};
\draw[fill=white] (0.6,0.5) circle (2pt) node[above]{$A$};
\draw[fill=white] (0.6,-0.5) circle (2pt) node[right]{$C$};
\draw[fill=white] (1.4,0.3) circle (2pt) node[below right]{$1$};
\draw[fill=white] (2.4,0.9) circle (2pt) node[below right]{$2$};
\draw[fill=white] (2.6,1.7) circle (2pt) node[above]{$4$};
\draw[fill=white] (1.6,1.1) circle (2pt) node[above]{$5$};
\draw[fill=white] (-0.5,0.8) circle (2pt) node[above]{$6$};
\draw[fill=white] (-1.3,0.6) circle (2pt) node[left]{$8$};
\draw[fill=white] (0,-1.1) circle (2pt) node[below left]{$9$};
\draw (0.3,-2.5) node {a};
\end{tikzpicture}\hspace{-1cm}
\begin{tikzpicture}
\tkzDefPoints{0/0/A,4/0/B,3/3/C};
\tkzDrawLine[very thick](A,B);
\tkzDrawLine[very thick](B,C);
\tkzDrawLine[very thick](C,A);
\tkzDefBarycentricPoint(B=1,C=2); \tkzGetPoint{L0};
\tkzDefBarycentricPoint(B=2,C=1); \tkzGetPoint{L9};
\tkzDefBarycentricPoint(C=1,A=5); \tkzGetPoint{L5};
\tkzDefBarycentricPoint(C=2,A=4); \tkzGetPoint{L4};
\tkzDefBarycentricPoint(C=3,A=3); \tkzGetPoint{L3};
\tkzDefBarycentricPoint(C=4,A=2); \tkzGetPoint{L2};
\tkzDefBarycentricPoint(C=5,A=1); \tkzGetPoint{L1};
\tkzDefBarycentricPoint(A=1,B=3); \tkzGetPoint{L8};
\tkzDefBarycentricPoint(A=2,B=2); \tkzGetPoint{L7};
\tkzDefBarycentricPoint(A=3,B=1); \tkzGetPoint{L6};
\tkzLabelPoint[above left=-2pt](L5){$5$};
\tkzLabelPoint[above left=-2pt](L4){$4$};
\tkzLabelPoint[above left=-2pt](L3){$3$};
\tkzLabelPoint[above left=-2pt](L2){$2$};
\tkzLabelPoint[above left=-2pt](L1){$1$};
\tkzLabelPoint[below](L8){$8$};
\tkzLabelPoint[below](L7){$7$};
\tkzLabelPoint[below right=-1pt](L6){$6$};
\tkzLabelPoint[above right=-2pt](L0){$10$};
\tkzLabelPoint[above right=-2pt](L9){$9$};
\tkzLabelPoint[below](A){$A$};
\tkzLabelPoint[below left=-2pt](B){$B$};
\tkzLabelPoint[right](C){$C$};
\tkzDrawLine(L2,L6);\tkzDrawLine(L5,L7);
\tkzInterLL(L2,L6)(L5,L7); \tkzGetPoint{X};
\tkzInterLL(L4,L0)(L9,L3); \tkzGetPoint{Y};
\tkzDrawLine(L0,L4); \tkzDrawLine(L9,L3);
\tkzLabelPoint[above right=-2pt](X){$X$};
\tkzLabelPoint[above](Y){$Y$};
\tkzDrawPoints[fill = white, size=3.5pt, ](A,B,C,L3,L4,L5,L6);
\tkzDrawPoints[fill = white, size=3.5pt, color=red ](L1,L2,L7,L8,L9,L0,X,Y);
\draw (2,-1) node {b};
\end{tikzpicture}\vspace{0.5cm}
\begin{tikzpicture}
\draw[->] (-0.25,0) -- (2.75,0);
\draw[->] (0,-0.25) -- (0,1.75);
\fill[color=gray!20](0,1.5)--(0,0)--(2.5,0)--cycle;
\draw[thick] (0,1.5)--(0,0)--(2.5,0);
\draw[thick,blue](0,1.5)--(2.5,0);
\foreach \x in {0,1,2,3,4,5} \fill (\x/2,0) circle(2pt);
\foreach \x in {0,1,2,3} \fill (\x/2,0.5) circle(2pt);
\foreach \x in {0,1} \fill (\x/2,1) circle(2pt);
\fill (0,1.5) circle(2pt);
\fill[white] (0,1) circle(1.5pt);
\fill[white] (0.5,1) circle(1.5pt);
\fill[white] (2,0) circle(1.5pt);
\draw (1.25,-0.8) node {c};
\end{tikzpicture}\qquad
\begin{tikzpicture}
\draw[->] (-1,0)--(1,0);
\draw[->] (0,-1)--(0,1);
\draw[-latex,blue] (-70:0.8)arc(-70:-50:0.8);
\draw[-latex,blue] (50:0.8)arc(50:70:0.8);
\draw[-latex,blue] (170:0.8)arc(170:190:0.8);
\fill[blue] (-70:0.8) circle(1.5pt);
\fill[blue]  (50:0.8) circle(1.5pt);
\fill[blue]  (170:0.8) circle(1.5pt);
\draw (-45:0.8) node {\tiny 8/3};
\draw (75:0.8) node {\tiny 8/3};
\draw (195:0.8) node {\tiny 8/3};
\draw (0,-1.5) node {d};
\end{tikzpicture}
\caption{}\label{e8}
\end{figure}

The versal deformation of the singularity $E_8$ is 
$$D(x,p)=x^5+p^3+a_1+a_2x+a_3x^2+a_4x^3+a_5p+a_6xp+a_7x^2p+a_8x^3p$$
the Newton polygon is shown of Fig.\ref{e8}c. It has one side $(-5,3)$ directed upward with multiplicity $1$. Thus the curve has genus 4 with one point at infinity with the rank of homology group 8. The growth diagram shown on Fig.\ref{e8}d consists of three points with angle $2\pi/3$ between them rotating with the angular speed $8/3$. It corresponds to the periodic word $[(s_1s_2)^8]$ and thus the Stokes data amounts to a configuration of $8$ points in $P^2$.

Exchanging $x$ and $p$ we get on the growth diagram $5$ points in a vertices of a regular pentagon rotating about its center with the angular speed $8/5$. It corresponds to a sequence of 8 points in $P^4$. The two configuration space are birationally isomorphic via Gale duality. 

\section{Other cases.}
We leave the detailed consideration of singularities of other types as an exercise. 

\begin{enumerate}
\item[$D_n$:]  The Dynkin diagram  corresponds to a configuration of triangles in the projective plane $P^2$ with 2, 2 and $n-2$ points on its respective sides. The corresponding singularity $xp^2+x^{n+1}$ corresponds to a configuration of $n$ lines in a two-dimensional local system on a circle. The differential operator corresponding to $p^{n+1}+x^2p$ corresponds to configurations of flags in $P^{n}$.
\item[$E_6$:] The Dynkin diagram as well as the singularity $p^3+x^4$ correspond to a configuration of triangles in $P^2$ with 2,3 and 3 points on the sides, respectively. The singularity $p^3+x^4$ corresponds to the word $[(s_1s_2)^7]$, i.e., to  configurations of 7-tuples of points of $P^2$. It is easy to see that the two configuration spaces are equivalent. The singularity $p^4+x^3$ corresponds to configurations of 7-tuples points in $P^3$.
\item[$E_7$:] The Dynkin diagram corresponds to a configuration of triangles in $P^2$ with 2,3 and 4 points on the sides, respectively. The singularity $p^3+x^3p$ corresponds to configurations of pentagons in $P^2$ with one marked point on each side, which is equivalent to the space of configurations given by the Dynkin diagram. The singularity $xp^3+x^3$ corresponds to a configuration space of 21-periodic sequences flags in $P^3$ which is too complicated to be described here. 
\end{enumerate}

\end{document}